\documentclass[a10paper,english,17pt]{article}
\usepackage[T1]{fontenc}
\usepackage{babel}
\usepackage{graphicx,color}
\usepackage{babel}
\usepackage{amsfonts}
\usepackage{amssymb, amsmath, amsfonts, amsthm, mathrsfs}
\usepackage{hhline, latexsym}
\usepackage{float}
\usepackage[margin=3cm]{geometry}
\begin{document}
\title{Blowing-up of locally monomially foliated space}
\author{Aymen Braghtha}
\maketitle
\begin{abstract}
In this paper, we prove that the blowing-up preserve the local monomiality of foliated space. 
\end{abstract}
\section{Locally monomial foliations}
Let $M$ be an analytic manifold  of dimension $n$ and $D\subset M$ be a divisor with normal crossings. We denote respectively by $\mathcal{O}_M$ and $\Theta_M[\log D]$ the sheaf of holomorphic functions and the sheaf of vector fields on $M$ which are tangent to $D$.

 A singular foliation on $(M,D)$ is coherent subsheaf $\mathcal{F}$ of $\Theta_M[\log D]$ which is reduced and integrable (see [1] and [2]). The dimension (or the rank) of $\mathcal{F}$ is given by 
$$
s=\max_{p\in M}\dim\mathcal{F}(p)
$$
where $\mathcal{F}(p)\subset T_pM$ denote the vector subspace  generated by the evaluation of $\mathcal{F}$ at $p$.

 Let $\mathbb{F}$ be a field (we usually take $\mathbb{F}=\mathbb{Q}, \mathbb{R}$ or $\mathbb{C}$). We shall say that $\mathcal{F}$ is \emph{$\mathbb{F}$-locally monomial} if for each point $p\in M$ there exists
\begin{enumerate}
\item a local system of coordinates $x=(x_1,\ldots,x_n)$ at $p$
\item an $s$-dimensional vector subspace $V\subset \mathbb{F}^n$ 
\end{enumerate}
such that $D$ is locally given by 
$$
D_p=\{x_i=0: i\in I\},\quad\text{for some}\quad I\subset\{1,\ldots,n\}
$$
and $\mathcal{F}_p$ is the $\mathcal{O}_{M,p}$-module generated by the abelian Lie algebra
$$
\mathcal{L}(V)=\{\sum^n_{i=1}a_ix_i\frac{\partial}{\partial x_i}: a\in V\}\bigoplus_{i\bar{I}:e_i\in M}\mathbb{F}\frac{\partial}{\partial x_i}
$$
where $\bar{I}=\{1,\ldots,n\}\setminus I$. We shall say that the triple $(M,D,\mathcal{F})$ is locally monomially foliated space and that $(x,I,V)$ is a local presentation at $p$.\\
\linebreak
\textbf{Lemma 1.1.} \emph{Let $(x,I,V)$ be a local presentation for $(M,D,\mathcal{F})$ at a point $p$. Then, for each vector $m\in V^{\perp}$, the (possibly multivalued) function $f(x)=x^m$ is a first integral of $\mathcal{F}$.}

\section{Blowing-up}
Let $Y\subset M$ be a smooth submanifold of codimension $r$. We shall say that $Y$ has \emph{normal crossings} with $(M,D,\mathcal{F})$ if for each point $p\in Y$ there exists a local presentation $(x,I,V)$ at $p$ such that $Y$ is given by 
$$
Y=\{x_1=x_2=\ldots=x_r=0\}.
$$
\textbf{Proposition 1.2.} \emph{Let $\Phi:\widetilde{M}\rightarrow M$ be the blowing-up with a center $Y$ which has normal crossings with $(M,D,\mathcal{F})$. Let 
$$
\widetilde{D}=\Phi^{-1}(D)\quad\text{and}\quad \widetilde{\mathcal{F}}\subset\Theta_{\widetilde{M}}[\log\widetilde{D}]
$$
denote the total transform of $D$ and the strict transform of $\mathcal{F}$ respectively. Then, the triple $(\widetilde{M},\widetilde{D},\widetilde{\mathcal{F}})$ is a locally monomially foliated space.
}\\

The proof is based on the following result on linear algebra.
\section{Some linear algebra}
Let $\mathbb{F}$ be a field and let $V\subset\mathbb{F}^n$ be a vector subspace of dimension $s$. Let us fix a disjoint partition of indices $\{1,\ldots,n\}=I_1\sqcup I_2$ write $\mathbb{F}^n=\mathbb{F}^{I_1}\oplus\mathbb{F}^{I_2}$ and let $\pi_{I}:\mathbb{F}^n\rightarrow\mathbb{F}^{I}$ denote the projection in the corresponding subspace $\mathbb{F}^I$ generated by $\{e_i:i\in I\}$.\\
\linebreak
\textbf{Lemma 3.1.} \emph{There exists a basis for $V$ such that for each vector $v$ in this basis, either
$$
v=\pi_{I_2}(v)\qquad\text{or}\qquad v=\pi_{I_2}(v)+e_i
$$
for some $i\in I_1$}
\begin{proof}
Up to a permutation of coordinates, we can suppose that $I_1=\{1,\ldots,n_1\}$ and $I_1=\{n_1+1,\ldots,n\}$ (with the convention that $n_1=0$ if $I_1=\emptyset$)

Let $M=[m_1,\ldots,m_s]$ be a $s\times n$ matrix whose rows are an arbitrary basis of $V$. By a finite number of elementary row operations and permutations of columns (which leave invariant the subsets $I_1$ and $I_2$), we can suppose that the matrix $M$ has the form
$$
M=
\left(\begin{array}{c|c||c|c} 
Id_{k_1\times k_1}&A_{k_1\times l_1}&0_{k_1\times k_2}&B_{k_1\times l_2}\\
\hline
0_{k_2\times k_1}&0_{k_2\times l_1}&Id_{k_2\times k_2}&C_{k_2\times l_2}
\end{array}\right)
$$
where $k_i+l_i=|I_i|$ for $i=1,2, k_1+k_2=s$ and $Id_{k,l}$ and $0_{k,l}$ denote the $k\times l$ identity and zero matrix respectively.

Now, it suffices to define $v_i=e_i+\sum_{n_1+k_2+1\leq j}m_{i,j}e_j$, for $i=1,\ldots,k_1$ and $v_i=e_{i+n_1}+\sum_{n_1+k_2+1\leq j}m_{i,j}e_j$ for $i=k_1+1,\ldots,s$.
\end{proof}

\textbf{Burgundy University, Burgundy Institue of Mathematics,\\
U.M.R. 5584 du C.N.R.S., B.P. 47870, 21078 Dijon \\
Cedex - France.}\\
\textbf{E-mail adress:} aymenbraghtha@yahoo.fr

\end{document}